\documentclass[bj]{imsart}

\usepackage{amsthm,amsmath,csquotes,amssymb}
\usepackage{floatrow,graphicx,tabu}
\newfloatcommand{capbtabbox}{table}[][\FBwidth]

\startlocaldefs
\newtheorem{thm}{Theorem}[section]
\newtheorem{cor}[thm]{Corollary}
\newtheorem{lem}[thm]{Lemma}

\theoremstyle{definition}
\newtheorem{defn}[thm]{Definition}

\newtheorem{remark}[thm]{Remark}
\renewcommand{\P}{\mathbb{P}}
\newcommand{\Q}{\mathbb{Q}}
\newcommand{\R}{\mathbb{R}}

\newcommand{\N}{\mathbb{N}}
\newcommand{\E}{\mathbb{E}}
\newcommand{\cG}{\mathcal{G}}
\newcommand{\cS}{\mathcal{S}}
\newcommand{\cH}{\mathcal{H}}

\endlocaldefs

\begin{document}

\begin{frontmatter}

\title{Degeneracy in sparse ERGMs with functions of degrees as sufficient statistics}

\runtitle{Degeneracy in sparse ERGMs}

\begin{aug}
 \author{\fnms{Sumit} \snm{Mukherjee}\thanksref{}\corref{}\ead[label=e1]{sm3949@columbia.edu}}

\address[]{Department of Statistics, Columbia University, New York, NY,
USA,\\ Research Partially Supported by NSF Grant DMS-1712037,\\
 {\printead{e1}}}

\runauthor{S. Mukherjee}

\affiliation{Columbia University}

\end{aug}

\begin{abstract}
A sufficient criterion for \enquote{non-degeneracy} is given for  Exponential Random Graph Models on sparse graphs with sufficient statistics which are functions of the degree sequence. This criterion explains why statistics such as alternating $k$-star are non-degenerate, whereas subgraph counts are degenerate. It is further shown that this criterion is \enquote{almost} tight. Existence of consistent estimates is then proved for non-degenerate Exponential Random Graph Models.
\end{abstract}





\end{frontmatter}

\section{Introduction}

Exponential families are frequently used in social science literature to model social networks  (see \cite{FS,H,HL,HH,MHH,SH,SPRH,WF,YRF}  and the references within).
Such models are usually referred to as Exponential Random Graph Models, commonly abbreviated  as ERGMs, in the social science community. Starting with   \cite{H}  in 2003, it has been noted in the social science literature  that ERGMs with subgraph counts do not behave in a nice manner in terms of sampling and estimation procedures. This phenomenon is typically referred to as degeneracy. Attempts have been made to characterize degeneracy (see for e.g. \cite{S,SR}) but there is no universally accepted definition for degeneracy.  This paper will adopt a notion similar to  \cite{H,SPRH}, where  the degeneracy of a model is attributed to the sufficient statistics of the model. That is,  the model will be deemed non-degenerate if the model behaves \enquote{nicely} for all choices of the parameter values. Thus under this notion, a degenerate model is caused by one or more degenerate statistics, and so  the term degenerate will be used for both the model as well as the statistic. 
\\

One of the features of degeneracy is that such models place most of their mass on a very small sub-collection of graphs. The intuitive idea behind their reasoning is that in such models the  neighboring edges are highly correlated. This causes a cascading effect through the graph, and so the model ends up putting most of its mass on very sparse or very dense graphs. In a sense, such models capture \enquote{too much interaction}. Thus an MCMC sample from such a model almost invariably gives either a very sparse graph, or a very dense graph. Another feature of such models is that small changes in the parameter can cause a large change in the underlying model. As such parameter estimates obtained from such models are usually not stable. 
\\

It has been subsequently noted in \cite{SPRH} in 2006 that not all ERGMs exhibit degeneracy in empirical studies. In fact, in this paper the authors argue that using modified versions of subgraph counts can reduce this problem to a large extent. The modifications are specifically aimed at reducing correlations between edges, and simulations seem to confirm this intuition. This raises the question of whether we can justify this empirical non degeneracy in a more rigorous setting, and whether we can develop an Inferential Framework for such models.

\subsection{Outline of the paper}

Section \ref{intro:eg} describes some examples of concrete interest and introduces the theoretical set up, and Section \ref{intro:main} outlines the main results of this paper.
Section \ref{examples} explains how one can use the results of this paper to compute normalizing constants using four examples. 

 The main tool for proving the results of this paper is a large deviation principle for the empirical degree distribution $\mu_n^G$  for a sparse Erd\"os-Renyi graph $G$ with respect to weak topology, studied in \cite[Corollary 2.2]{DM} and  \cite[Corollary 1.9]{BC}.  Their result is outlined in Section \ref{ldp}. Section \ref{sec:proofs}  carries out the proofs of the main results of the paper (Theorem \ref{thm:infty}, Corollary \ref{gwd}, and Theorem  \ref{gwd2}), using auxiliary lemmas which are proved in Section \ref{sub:proofs}. Finally, Section \ref{fit_the_model} proves existence of consistence estimates for the ERGMs proposed in this paper (c.f. Theorem \ref{showconsistency1} and \ref{showconsistency2}).

\subsection{Examples}\label{intro:eg}

The following definition gives the necessary notations for introducing some of the examples from \cite{SPRH} which are non degenerate at an empirical level.

\begin{defn}

Let $\mathcal{G}_n$ denote the space of all simple labelled undirected graphs on $n$  vertices. 
For any $G\in \mathcal{G}_n$ let ${\bf d}(G)=(d_1(G),\cdots,d_n(G))$ denote the labeled degree sequence of $G$, i.e. $d_j(G)$ is the degree of vertex $j$. Also let $E(G):=\frac{1}{2}\sum\limits_{j=1}^nd_j(G)$ denote the number of edges in $G$.

For $0\le i\le n-1$, let $h_i(G):=\#\{1\le j\le n:d_j(G)=i\}$ denote the number of vertices of degree $i$. Summing over $i$ gives  $\sum\limits_{i=0}^{n-1}h_i(G)=n$, since the sum is over all the vertices of $G$.  The quantity ${\bf h}(G):=\{h_i(G)\}_{i=0}^{n-1}$  will be referred to as the degree frequency vector.

Recall that a $k$-star has $k$ edges and $k+1$ vertices. For any $k\ge 2$, let $T_k(G)$ denote the number of copies of $k$-stars in $G$. The counting scheme is such that all copies of the $k$-star are considered, and not just the induced ones. This counting scheme gives the following simple formula for $T_k(G)$ in terms of its degrees ${\bf d}(G)$, as well as the degree frequency vector ${\bf h}(G)$:
$$T_k(G)=\sum\limits_{j=1}^n{d_j(G)\choose k}=\sum\limits_{i=0}^{n-1}h_i(G){i\choose k}.$$
This is because for any vertex $j$, there are ${d_j(G)\choose k}$ $k$-stars with $j$ as the center vertex, and so  adding over $j$ gives the total number of $k$-stars. The second equality follows by rearranging the first sum.

\end{defn}

We will now introduce some of the non-degenerate statistics defined in \cite{SPRH}.
\begin{enumerate}
\item[(a)]{\bf Geometrically weighted degree statistic}

The geometrically weighted degree statistic has the form $$\text{gwd}_\alpha(G):=\sum\limits_{i=0}^{n-1}e^{-\alpha i}h_i(G),$$ where $\alpha>0$ is known. The geometrically decaying weights ensure that  the contribution of vertices with large degree is negligible. Thus as the degrees of the graph increase, the statistic does not grow too fast, and cascading effect of this statistic is reduced.
\\

\item[(b)]{ \bf The alternating $k$-star}

For a fixed parameter $\lambda>1$, the alternating $k$-star is defined as
$$\text{aks}_\lambda(G):=\sum\limits_{k=2}^{n-1}\frac{(-1)^k}{\lambda^{k-2}}T_k(G),$$
where $T_k$'s are the $k$-star counts defined above. 
In this case again the geometrically decaying weights ensure that the cascading effects of higher star counts is reduced. Also because of the alternate signs the cascading effect of consecutive terms is cancelled to a large extent.

The authors in \cite{SPRH} note that using the formula for $T_k(G)$ in terms of the degree frequency vector ${\bf h}(G)$, the alternating $k$-star statistic can be written as
$$\text{aks}_\lambda(G)=\lambda^2\sum\limits_{i=0}^{n-1}\Big[\Big(1-\frac{1}{\lambda}\Big)^i-1+\frac{i}{\lambda}\Big]h_i(G)=\lambda^2\text{gwd}_\alpha(G)-n\lambda^2+2\lambda E(G)$$
with $e^{-\alpha}=1-1/\lambda$. Thus the two statistics $\text{gwd}_\alpha$ and $\text{aks}_\lambda$ are connected by a simple formula, and both these statistics are functions of the degree frequency vector ${\bf h}$. 
\\

\item[(c)]{\bf The number of isolated nodes}

The statistic $h_0(G)$ which is the number of isolated vertices in the graph $G$. This statistic is obtained from the $\text{gwd}_\alpha$ statistic by letting $\alpha\rightarrow \infty$, or equivalently from the $\text{aks}_\lambda$ by letting $\lambda\rightarrow 1$. 
\\

\item[(d)]{{\bf The Yule distribution statistic}}

Another statistic which penalizes high degrees is the Yule distribution statistic, given by
$$yu(G):=\sum_{j=1}^n\frac{1}{(d_j+c)_r}=\sum_{i=0}^{n-1}\frac{1}{(i+c)_r}h_i(G),\quad (d)_r:=d(d+1)\cdots(d+r-1),$$
where $r$ and $c$ are both positive integers. In this case the penalty is polynomial as opposed to geometric as in $\text{gwd}_\alpha$, but a similar non-degenerative effect is achieved.
\\

\end{enumerate}

In all the four examples above the statistic under consideration can be written as $\sum_{i=0}^{n-1}f(i)h_i(G)$ for some function $f:\mathbb{N}_0\mapsto \mathbb{R}$, where $\mathbb{N}_0:=\mathbb{N}\cup \{0\}$.  As an illustration, the $\text{gwd}_\alpha, \text{aks}_\lambda$, the number of isolated vertices and the Yule distribution statistic fit this framework  with $$f(i)= e^{-\alpha i},\quad f(i)=\lambda^2\Big[\Big(1-\frac{1}{\lambda}\Big)^i-1+\frac{i}{\lambda}\Big],\quad f(i)=1_{i=0},\quad f(i)=\frac{1}{(i+c)_r}$$ respectively. Restricting attention to statistics of this form, one can ask when is this statistic well behaved. The results of  this paper gives a sufficient condition for this which is easy to check:
\begin{align}\label{consistency}
\lim_{i\rightarrow\infty}\frac{|f(i)|}{i\log i}=0.
\end{align}
This will be made precise in Theorem \ref{thm:infty} and Corollary \ref{gwd}.

In particular, \eqref{consistency} holds for $\text{gwd}_\alpha$ for $\alpha>0$, and the number of isolated vertices, and the Yule distribution statistic, as in all these cases the function $f$ is bounded. For $\text{aks}_\lambda$ with $\lambda>1$ the function  $$f(i)=\lambda^2\Big[\Big(1-\frac{1}{\lambda}\Big)^i-1+\frac{i}{\lambda}\Big]$$ is unbounded, but dominated by the linear term. 
On the other hand, the number of $k$-stars also equals $\sum_{i=0}^{n-1}f(i)h_i(G)$ for the choice $f(i)={i\choose k}$ which does not satisfy \eqref{consistency}, as in this case $f(i)$ grows at a polynomial rate. Also in the  alternating $k$-star statistic if the signs do not alternate, then one has
$$\sum\limits_{k=2}^{n-1}\frac{1}{\lambda^{k-2}}T_k(G)=\sum_{i=0}^{n-1}f(i)h_i(G)$$
with $f(i)=\lambda^2\Big[\Big(1+\frac{1}{\lambda}\Big)^i-1-\frac{i}{\lambda}\Big]$ which does not satisfy \eqref{consistency}, as in this case the exponential term dominates. Thus it is crucial that the signs in the alternating $k$-star statistic do alternate. It follows from Theorem \ref{gwd2} that both the number of $k$-stars and the non-alternating $k$-star statistics are degenerate, in a sense which is again made precise in Theorem \ref{gwd2}.

 The main tool for these results is the analysis of sparse graphs, as opposed to dense graphs as in  \cite{BBS,CD,RY}. Recall that  in a dense graph on $n$ vertices the number of edges is $O(n^2)$ and  the degrees are $O(n)$.  Here and henceforth in this paper, we use the notation $a_n=O(b_n)$ for two positive real sequences $\{a_n\}_{n\ge 1}$ and $\{b_n\}_{n\ge 1}$, if there exists a constant $C$ free of $n$ such that $a_n\le C b_n$ for all $n\ge 1$. With this notation, the term sparse graphs will refer to graphs which have $O(n)$ edges and the degrees of the vertices are $O(1)$. One reason it is interesting to model sparse graphs is that most real life networks seem to be sparse.  Another reason is that the dense graph theory does not provide a good explanation for why the modified versions of subgraph counts mentioned above (such as $\text{aks}_\lambda$) are non-degenerate, whereas the subgraph count statistics (such as star counts) are.

For a unified treatment of these and other examples, consider an exponential family on $\mathcal{G}_n$ of the form 
\begin{align}\label{thesis0}
\Q_{n,\beta,f}(G):=\Big(\frac{\beta}{n}\Big)^{E(G)}\Big(1-\frac{\beta}{n}\Big)^{{n\choose 2}-E(G)}e^{ \sum\limits_{i=0}^{n-1}h_i(G)f(i)-Z_n(\beta,f)},
\end{align}
where $f:\mathbb{N}_0\mapsto \R$, $\beta$ is a positive real valued parameter, and $Z_n(\beta,f)$ is the log normalizing constant.
If  $f$ is either identically $0$ or exactly linear, this model reduces to a sparse  Erd\"os-Renyi model  which puts most of its mass on sparse graphs. Thus the same should be true for functions $f(.)$ which do not grow too fast. Since the model does not change if $f(.)$ is replaced by $f(.)+c$ for some constant $c$, without loss of generality we will assume $f(0)=0$.

It should be noted at this point that $\Q_{n,\beta,f}$ is not the same as the ${\bf \beta}$ model studied in  \cite{CDS}. The ${\bf \beta}$ model is an exponential family on  $\mathcal{G}_n$ whose probability mass function is proportional to $\text{exp}\{\sum\limits_{j=1}^n\beta_jd_j(G)\}$ where ${\bf \beta}=(\beta_1,\cdots,\beta_n)$ is an $n$ dimensional parameter.  In the $\beta$ model the labeled degree sequence $(d_1(G),\cdots,d_n(G))$ is minimal sufficient. On the other hand in the model of  \eqref{thesis0} if the function $f$ is assumed to be unknown the minimal sufficient statistics are the unlabeled  degree sequence $(d_{(1)}(G)\ge d_{(2)}(G)\ge \cdots \ge d_{(n)}(G))$, or equivalently the degree frequency vector $(h_0(G),h_1(G),\cdots,h_{n-1}(G))$. More importantly, model \eqref{thesis0} introduces non trivial dependence among the edges of the graph $G$, whereas under the ${\bf \beta}$ model the edges are mutually independent. 
In \cite{CDS} the authors worked in the dense graph regime and showed that if the components of the parameter vector ${\bf \beta}$ stays uniformly bounded, then all entries of ${\bf \beta}$ can be simultaneously estimated consistently. In a similar manner, Theorem \ref{showconsistency2} shows  that if the true function $f$ is unknown and treated as a parameter, one can estimate the value of the function $f$ consistently at every fixed $i$, under the assumption that $f$ satisfies  \eqref{consistency}.


\subsection{Statement of Main results}\label{intro:main}
For analyzing model \eqref{thesis0} it  suffices to study the degree sequence. The following definition encodes the entire degree sequence as one probability measure on non-negative integers.

\begin{defn}
Given the labelled degrees of a graph $(d_1(G),\cdots,d_n(G))$, the empirical distribution of the degree sequence is defined by $\mu_n^G:=\frac{1}{n}\sum\limits_{j=1}^n\delta_{d_j(G)}$ i.e. $\mu_n^{G}$ is the measure which puts mass $1/n$ at each of the observed degree $d_j(G)$, and  is a probability measure on $\mathbb{N}_0$.
An equivalent definition of $\mu_n^{G}$ in terms of the degree frequency vector ${\bf h}(G)$ is the probability measure which puts mass $h_i(G)/n$ at $i$, for $0\le i\le n-1$.

 With this definition, any statistic of the form $\sum\limits_{i=0}^{n-1}f(i)h_i(G)$ can be written as $n\mu_n^{G}[f]$, where $\mu[f]$ denotes the mean of $f$ with respect to the measure $\mu$ (when it exists), i.e.
$$\mu[f]:=\sum\limits_{i=0}^\infty \mu(i)f(i).$$In particular if $f(i)=i$ is the identity function, then define
$\overline{\mu}:=\sum_{i=1}^\infty i\mu(i)$ to be the mean of the measure $\mu$.

\end{defn}

The next definition gives all the necessary ingredients for expressing the asymptotic log normalizing constant as a one dimensional optimization problem.

\begin{defn}\label{dd}
Suppose the function $f:\N_0\mapsto \R$  satisfies 
\begin{align}\label{eq:consistency2}
\limsup_{i\rightarrow\infty}\frac{f(i)}{i\log i}=0,
\end{align}
which is a slightly weaker condition than \eqref{consistency}. For $u\ge 0$ define an exponential family on $\mathbb{N}_0$ with probability mass function
$$\sigma_{u,f}(i)=\frac{1}{i!}u^i e^{f(i)-Z(u,f)},$$
where $Z(u,f)$ is the log normalizing constant, i.e. 
$$Z(u,f):=\log\Big(\sum\limits_{i=0}^\infty \frac{1}{i!}u^i e^{f(i)}\Big).$$ Since $f$ satisfies \eqref{eq:consistency2} we have that $Z(u,f)<\infty$. 
 Let $\Omega_{f}$ denote the set of all probability measures of the form $\sigma_{u,f}$ for $u\ge 0$. Also let $m(u,f):=\overline{\sigma_{u,f}}$ denote the mean of $\sigma_{u,f}$.
 Finally, for $\beta>0$ let $J(\beta,f)$ denote the solution to the following optimization problem  
\begin{align}\label{eq:j}
J(\beta,f):=\sup_{u  \ge 0}\Big\{Z(u,f)-m(u,f)\log u +\frac{m(u,f)}{2}\log (m(u,f)\beta)-\frac{m(u,f)+\beta}{2}\Big\}.
\end{align}
The definition of $J(\beta,f)$ involves an optimization  over the scalar non-negative variable $u$ which can be computed numerically.

\end{defn}

The first main result of this paper is the following theorem, which gives the asymptotics of the log normalizing constant for the model $\Q_{n,\beta,f}$ under assumptions on the growth rate of $f$. Existence of limiting log normalizing constant for a dependent system with growing number of variables governed by a Gibbs measure has attained considerable interest in Statistical Physics, where this is typically referred to as existence of the thermodynamic limit. Typically the limiting normalizing constant is expressed in terms of an optimization problem, and the optimizers represent the steady states of the distribution. See \cite{Ruelle} for more on existence of thermodynamic limits and its properties in general.

\begin{thm}\label{thm:infty}

Suppose either of these two conditions hold:
\begin{enumerate}
\item[(i)] 
$f:\mathbb{N}_0\mapsto\R$ satisfies \eqref{consistency}, 

or

\item[(ii)]
$f:\mathbb{N}_0\mapsto\R$  is non increasing.

\end{enumerate}

Let $G$ be a random graph  from the exponential family $\Q_{n,\beta, f}$   as defined in \eqref{thesis0}.

 \begin{enumerate}
 \item[(a)] 
 Then as $n\rightarrow\infty$, the asymptotics of the log normalizing constant is given by

$$\lim\limits_{n\rightarrow\infty}\frac{1}{n} Z_n(\beta,f)=J(\beta, f),$$
with $J(\beta,f)$ as in definition \ref{dd}.

\item[(b)]
   The supremum in the definition of $J(\beta, f)$ is attained on a finite set of positive reals $\{u_1,u_2\cdots,u_k\}$  satisfying the equation $u_l^2=\beta \bar{\sigma}_{u_l, f}$ for $1\le l\le k$, with $\sigma_{u,f}$ as in definition \ref{dd}. Further, for any function $\psi$ satisfying \eqref{consistency} one has  $$\min_{i=1}^k \Big|\mu^G_n(\psi)-\sigma_{u_i, f}(\psi)\Big|\stackrel{p}{\rightarrow}0.$$
where $\mu_n^{G}$ is the empirical degree distribution of $G$.
\end{enumerate}
\end{thm}

Note that both the conditions  (i) and (ii) considered in part (a) of Theorem \ref{thm:infty} are sub-cases of the assumption \eqref{eq:consistency2}, which is used to ensure that $J(\beta,f)$ introduced in \eqref{eq:j} is well defined and finite. It is possible that the conclusion of Theorem \ref{thm:infty} holds for all $f$ satisfying \eqref{eq:consistency2}. 

An immediate application of the above theorem gives the following corollary.

\begin{cor}\label{gwd}
 
Suppose $f:\mathbb{N}_0\mapsto \R$ satisfy \eqref{consistency}, and
let $G$ be a random graph  from the exponential family $\Q_{n,\beta,\theta f}$, where $\Q_{n,\beta,f}$ is   as defined in \eqref{thesis0}. Then the following conclusions hold:
\begin{enumerate}
\item[(a)]
Both part (a) and part (b) of Theorem \ref{thm:infty} hold with $f$ replaced by $\theta f$, for all $\theta\in \R$. Also, the limiting log partition function 
 $$J(\beta,\theta f)=\lim\limits_{n\rightarrow\infty}\frac{1}{n} Z_n(\beta,\theta f)$$
 is finite and continuous in $\theta$. 
 
 \item[(b)]
 There exists positive reals $m<M$ depending on $(f,\beta,\theta)$ such that $$\lim_{n\rightarrow\infty}\Q_{n,\beta,\theta f}\Big(\frac{E(G)}{n}\in[m,M]\Big)=1.$$
 \end{enumerate}

\end{cor}

\begin{remark}\label{lmt}
Part (a) of Corollary \ref{gwd} says that if  $|f|$ grows at a rate smaller than $i\log i$, then the corresponding model $\Q_{n,\beta,\theta f}$ is well behaved for both positive and negative $\theta$, in the sense that the limiting log partition function is finite and continuous in $\theta$. It also shows that the empirical degree distribution $\mu_n^{G}$ roughly behaves like a mixture of $\{\sigma_{u_i,\theta f}\}_{i=1}^k$ for large $n$. In particular if there is a unique optimizer $u_0$ to the optimization problem $J(\beta,\theta f)$, then the empirical degree distribution $\mu_n^G$ converges weakly to $\sigma_{u_0,\theta f}$, and $\bar{\mu}_n^G$ converges to $\bar{\sigma}_{u_0,\theta f}$.
\\

Part (b) shows that irrespective of whether there is a \enquote{phase transition}, the number of edges is linear in the number of vertices for all parameter values $\theta$ (c.f. \cite{Ruelle} for details on phase transitions in models of Statistical Mechanics). Thus the level of sparsity of the graph does not change with the parameter.
\\

Also, none of the limit points of the degree distribution is a Poisson, as $\sigma_{u,\theta f}$ is not a Poisson distribution unless $f$ is identically $0$ or linear, in which case the model $\Q_{n,\beta,\theta f}$ itself is a sparse Erd\"os-Renyi graph. On the other hand, the empirical degree distribution of a sparse Erd\"os-Renyi graph converges to Poisson. Thus unlike ERGMs on dense graphs as studied in \cite{CD}, ERGMs on sparse graphs do not behave like mixture of Erd\"os-Renyi graphs. Also, in the case of sparse ERGMs, it is possible to estimate multiple parameters consistently from a large single graph. In particular, see Theorem \ref{showconsistency1}  which constructs consistent estimates for $(\beta,\theta)$ when $f$ is known, and Theorem \ref{showconsistency2} which constructs consistent estimates for the function $f$ if $f$ is unknown. It should be noted here that consistent estimation of parameters in ERGMs was achieved in \cite{SR}, but under the assumption that the ERGM restricted to $n$ vertices is a projection of the corresponding ERGM on $n+1$ vertices. Consistency results have also been obtained for sparse ERGMs in \cite{KK}, but here the authors assume dyadic independence. In contrast, the models presented in this paper are neither projective nor have dyadic independence, and yet consistent estimation is possible in this case.
\\

 Since choosing a function $f$ is equivalent in spirit to specifying the degree distribution of the graph, one can fit a wide class of degree distributions by choosing a corresponding function $f$. Of course restriction \eqref{consistency} ensures that the degree distribution will have a finite exponential moment, which rules out degree distribution with power law tails. Power law tails correspond to the case when $f(i)$ grows at the rate $i\log i$, which require a more delicate analysis and is not carried out in this paper. 
\end{remark}

The next Theorem shows that some growth condition on $f$ needs to satisfied for the model $\Q_{n,\beta,\theta f}$ to be well behaved for all values of $\theta$. 
 \begin{thm}\label{gwd2}
Suppose $f:\mathbb{N}_0\mapsto \R$ is a non-decreasing function, and $G$ be a random graph from  the exponential family  $\Q_{n,\beta,\theta f}(.)$, where $\Q_{n,\beta,f}$ is as defined in \eqref{thesis0}.

\begin{enumerate}
\item[(a)]
If $\theta< 0$ then both parts (a) and part (b) of Theorem \ref{thm:infty} hold with $f$ replaced by $\theta f$. Also, the asymptotic log normalizing constant
 $$J(\beta,\theta f)=\lim\limits_{n\rightarrow\infty}\frac{1}{n} Z_n(\beta,\theta f)$$
 is finite and continuous in $\theta$. Further, there exists positive constants $m<M$ depending on $(\theta,\beta,f)$ such that $$\lim_{n\rightarrow\infty}\Q_{n,\beta,\theta f}\Big( \frac{E(G)}{n}\in [m,M]\Big)=1.$$
\\

\item[(b)]
If $f$ further satisfies \begin{align}\label{eq:dir}
\liminf_{i\rightarrow\infty}\frac{f(i)-f(i-1)}{\log i}>4,
\end{align} then for $\theta>0$ we have
$$\lim\limits_{n\rightarrow\infty}\frac{1}{nf(n)}Z_n(\beta,\theta f)=\theta,$$
and $\lim_{n\rightarrow\infty}\Q_{n,\beta,\theta f}(G=K_n)=1$.

\end{enumerate}
\end{thm}

\begin{remark}\label{degen}
Note that the assumption \eqref{eq:dir} automatically implies $f(i)$ is at least of order $i\log i$, i.e. \eqref{consistency} does not hold. Under this assumption, Theorem \ref{gwd2} demonstrates degeneracy in the sense of \cite{H,SPRH} in two ways. First, in this case the behavior of the model $Q_{n,\beta,\theta f}$ changes drastically at the origin.  For $\theta< 0$  the model puts all its mass on sparse graphs with $O(n)$ edges,  whereas for $\theta>0$ the model suddenly shifts all its mass  to the complete graph where number of edges is ${n\choose 2}\sim \frac{n^2}{2}$. Also, for $\theta>0$ the model puts most of its mass on a very small subset of $\cG_n$ (namely a subset of size $1$). Thus model \eqref{thesis0} can indeed be degenerate without any growth conditions on $f$.

In particular this happens for the choice $f(i)={i\choose k}$ for any $k$ fixed, for which the statistic $\sum_{i=0}^{n-1}f(i)h_i(G)$ becomes the number of $k$-stars, and for the choice $f(i)=\lambda^2\Big[\Big(1+\frac{1}{\lambda}\Big)^i-1-\frac{i}{\lambda}\Big]$, for which the statistic $\sum_{i=0}^{n-1}f(i)h_i(G)$ is the non alternating $k$-star. Note that in both these cases, the function $f$ is indeed non-decreasing, and satisfies \eqref{eq:dir}.


 \end{remark}

\subsection{Identifiability and estimating  parameters}

Since model $\Q_{n,\beta, \theta f}$ is well behaved for all $\theta$ when $f$ satisfies \eqref{consistency}, this subsection explores the estimation of parameters of the model, under the assumption that $f$ satisfies \eqref{consistency}. 
Assuming that $f$ is known, one can focus on estimating  the parameters $(\beta,\theta)$ in the model $\Q_{n,\beta,\theta f}$
from one sample $G$ from this model. 
If $f$ is exactly linear, i.e. there exists  a constant $b$ such that $f(i)=bi$, then the model $\Q_{n,\beta,\theta f}$ is same as Erd\"os-Renyi with parameter $$\frac{1}{1+\frac{n-\beta}{\beta}e^{-\theta b} }\approx \frac{\beta e^{\theta b}}{n}.$$
This model is asymptotically not identifiable along the curve where $\beta e^{\theta b}$ is constant, and so joint estimation of both parameters $(\beta,\theta)$ is not possible. If $f$ is not linear, consistent estimation  of both the parameters is possible under this model.

In order to motivate our proposed estimates, recall the prediction of part (b) of Corollary \ref{gwd},  that for large $n$ we have  $$\frac{h_i(G)}{n}\approx \frac{u^i}{i!}e^{\theta f(i)-Z(u,\theta f)},$$ if $f$ satisfies \eqref{consistency}. Taking this to be an exact equality, multiplying both sides by $i!$ and taking $\log$ gives
$$\log \frac{i! h_i(G)}{n}=- Z(u,\theta f)+\theta f(i)+i\log u.$$
Thus taking $x_1(i)=f(i)$, $x_2(i)=i$, $y(i)=\log \frac{i! h_i(G)}{n}$, we get a linear equation of the form $$y(i)=-Z(u,\theta f)+\theta x_1(i)+(\log u) x_2(i),$$
and so by fitting a multiple linear regression model using least squares with $y$ as response and $\{x_1,x_2\}$ as explanatory variables we can estimate $\theta$ and $\log u$. Finally note that $u,\theta,\beta$ are connected by the equation $u^2=\beta \bar{\sigma}_{u,\theta f}$, as shown in part (b) of Corollary \ref{gwd}. Since the empirical average of degrees $\bar{d}(G)=\frac{2E(G)}{n}$ converge to $\bar{\sigma}_{u,\theta f}$ in probability, one can use the approximate equation $nu^2=2\beta  E(G)$ along with the least squares estimate of $\log u$ to get an estimate of $\beta$. We will now show that the estimates of $(\theta,\beta)$ outlined above are indeed consistent.

\begin{thm}\label{showconsistency1}
Let $f:\mathbb{N}_0\mapsto\R$ be a known function which satisfies \eqref{consistency}, and $\theta_0\in \R, \beta_0>0$ be the true unknown parameters. Let $L$ be a fixed positive integer free of $n$ such that $f(i)/i$ is not constant for all $i\in [0,L]$.
Let $(\hat{\theta}_n,\hat{u}_n)$ be the least square estimates of $(\theta,u)$ defined via the following optimization problem:
\begin{align*}
(\hat{\theta}_n,\hat{u}_n):=\arg\inf_{c,\theta\in \R,u>0}\sum_{i=0}^L\Big\{\log\frac{i!h_{i}(G)}{n}-c-\theta f(i)-i\log u\Big\}^2.
\end{align*}
Then as $n\rightarrow\infty$, one has $\hat{\theta}_n\stackrel{p}{\rightarrow}\theta_0$. Further, the estimator $\hat{\beta}_n:=\frac{n\hat{u}_n^2}{2E(G)}\stackrel{p}{\rightarrow}\beta_0$.
\end{thm}

\begin{remark}\label{rem:new}
The estimates $(\hat{\theta},\hat{\beta})$ of the previous theorem are motivated by the fact that $$\frac{h_i(G)}{n}\approx \frac{1}{e^{Z(u,\theta f)}} e^{\theta f(i)} u^i,$$ as predicted in Corollary \ref{gwd} under the assumption that $f$ satisfies \eqref{consistency}.

Even though one can use a larger value of $L$ (for e.g. $L=n-1$), estimates of $f(i)$ for large $i$ are not as reliable. In particular, for large $i$ it is possible to have $h_i(G)=0$ which will give undefined values for $(\hat{\theta}_n,\hat{u}_n)$. This is the reason for choosing $L$ fixed, free of $n$. Given a graph $G$, any valid choice of $L$ must satisfy  $L\le L_n(G):=\max_{1\le j\le n}d_j(G)$. Indeed this is because $h_i(G)=0$ for $i>L_n(G)$, and so the least square optimization problem in Theorem \ref{showconsistency1} is not defined. A natural choice of $L$ is the maximum $i$ such that $h_i(G)$ is non zero for all $i\in [0,L]$.

Frequently it is the case that an observed graph $G$ has no isolated vertices. For example any person in a social network has at least one friend. A natural model in this case is the same exponential family, but conditioned to have no isolated vertices. Since $h_0(G)=0$, the estimator defined in Theorem \ref{showconsistency1} becomes undefined. In such cases, instead of starting at $0$ one can consider the values $i\in [1,L]$ for the least squares procedure. The same proof shows that the resulting estimator is consistent, whenever $f$ satisfies \eqref{consistency}.
\end{remark}

If there is no reasonable guess for the function $f$, then one can think of estimating the whole function $f$ in the model $\Q_{n,\beta, f}$ using one large graph $G$.  For any $f$ the two models 
$\Q_{n,\beta, f}$ and $\Q_{n,1,\tilde{f}}$ are asymptotically unidentifiable, where $$\tilde{f}(i)= f(i)+(i/2)\log \beta,$$
and so without loss of generality one may further assume $\beta=1$.  Under these assumptions, the next theorem reconstructs the whole function $f$.

\begin{thm}\label{showconsistency2}
Let $f:\mathbb{N}_0\mapsto \R$ be such that $f$ satisfy \eqref{consistency}, and consider the model $\Q_{n,1,f}$ where $\Q_{n,\beta,f}$ is as defined  in \eqref{thesis0}. Setting  $\hat{u}_n:=\sqrt{\frac{2E(G)}{n}}$  the function $\hat{f}_n:\mathbb{N}_0\mapsto\R$ defined by
$$\hat{f}_n(i)=\log\Big[\frac{i!h_i(G)}{h_0(G)}\Big]-i\log \hat{u}_n$$
satisfies 
$$\hat{f}_n(i)\stackrel{p}{\rightarrow}f(i)$$ for $i\ge 1$, as $n\rightarrow\infty$.

\end{thm}

\subsection{Scope for future work}

The statistics considered in this paper are functions of the degree sequence, or equivalently functions of $1$ neighborhoods of the graph. The literature has also focused on statistics which cannot be expressed in terms of the degrees, for example the alternating $k$-triangle statistic (for more details on this statistic refer to \cite{H,HH,SPRH}). The alternating $k$-triangle statistic depends on $2$ neighborhoods of a vertex and not $1$. The  large deviation result of \cite{BC} applies for any finite neighborhood, and so it seems plausible that the two neighborhood can be dealt with a modified version of the strategy of this paper. Of course, for $2$ (and general) neighborhoods the involved rate function will be more complicated.

\section{Some examples}\label{examples}

This Section uses the results of this paper to analyze four ERGMs on sparse graphs from the probability mass function $\Q_{n,\beta,\theta f}$ of  \eqref{thesis0}.  To specify the model it suffices to choose the function $f$. Note that in none of these models a closed form expression for the normalizing constant $Z_n(\beta,\theta f)$ seems available. Using Corollary \ref{gwd}, one can get numerical approximations for the asymptotic normalizing constant.

\subsection{Geometrically weighted degree}
In this case we have $\text{gwd}_\alpha=\sum_{i=0}^{n-1}h_i(G)f(i)$ with  $f(i)=e^{-\alpha i}$ for some $\alpha>0$. An application of Corollary \ref{gwd} gives the asymptotics of the log normalizing constant as
$$\lim_{n\rightarrow\infty}\frac{1}{n}Z_n(\beta, \theta f)=\sup_{u\ge 0}\Big\{Z(u,\theta f)-m(u,\theta f)\log u+\frac{m(u,\theta f)}{2}\log(m(u,\theta f)\beta)-\frac{m(u,\theta f)+\beta}{2}\Big\},$$
where $Z(u, \theta f)$ and $m(u, \theta f)$ are the log normalizing constant and mean respectively, of the probability mass function $\sigma_{u,\theta  f}$ on non-negative integers given by
$$\sigma_{u, \theta f}(i)\propto \frac{1}{i!}u ^i e^{\theta  f(i)}.$$
Setting $\gamma:=e^{-\alpha}$ one has
\begin{align*}
e^{Z(u,\theta f)}=\sum_{i=0}^\infty \frac{u^i}{i!}\text{exp}\{ \theta\gamma^i\}=\sum_{i=0}^\infty \frac{u^i}{i!}\sum_{j=0}^\infty \frac{\theta^j\gamma ^{ij}}{j!}=e^{u+\theta}\sum_{i,j=0}^\infty \frac{e^{-u}u^i}{i!} \frac{e^{-\theta}\theta^j}{j!}\gamma^{ij}=
e^{u+\theta}\E \gamma^{ XY},
\end{align*}
where $X,Y$ are mutually independent and $X\sim Pois(u),Y\sim Pois(\theta)$. By a similar calculation one has
\begin{align*}
m(u,\theta f)=\frac{\sum_{i=1}^\infty i\frac{u^i}{i!}\text{exp}\{\theta\gamma^i\}}{\sum_{i=0}^\infty \frac{u^i}{i!}\text{exp}\{ \theta\gamma^ i\}}
=u\frac{\sum_{i=0}^\infty \frac{u^i}{i!}\text{exp}\{ \theta\gamma  \gamma^{i}\}}{\sum_{i=0}^\infty \frac{u^i}{i!}\text{exp}\{\theta \gamma^i\}}
=u \frac{ e^{u+\theta\gamma }E \gamma^{ XZ}}{e^{u+\theta}E \gamma^{XY}}=u  e^{\theta(\gamma-1)}\frac{\E \gamma^{ XZ}}{\E \gamma^{XY}},
\end{align*}
where $Z\sim Poisson(\theta \gamma)$ independent of $X$.

Since closed form expressions are not known for moment generating function of products of independent Poissons, further simplification is not possible in this case. Of course one can use numerical approximations  by simulating an i.i.d.  sample of products of Poissons, and then using strong law of large numbers to estimate the moment generating function.

 \subsection{Logarithmic model}
 For this model set $$f(i)=-\log (i+1)_r=-\log(i+1)(i+2)\cdots(i+r),$$ where $r$ is a positive integer. In this case $|f(i)|$ grows logarithmically, and so by Corollary \ref{gwd} the asymptotics of the normalizing constant requires only the knowledge of $Z(u,\theta f)$ and $m(u,\theta f)$. For the special case $\theta=1$, a direct computation shows that 
 $$e^{Z(u,f)}=\frac{1}{u^r}\Big[e^u-\sum_{i=0}^{r-1}\frac{u^i}{i!}\Big],\quad m(u,f)=u-r+\frac{\frac{u^{r}}{(r-1)!}}{e^u-\sum_{i=0}^{r-1}\frac{u^i}{i!}}$$
 Thus in this case both $Z(u,f)$ and $m(u,f)$ are explicit, and numerical optimization of $J(\beta,f)$ is easy to carry out. No such simple formula exists for $Z(u,\theta f)$ and $m(u,\theta f)$ for $\theta\ne 1$.
 
\subsection{Sparse penalty model}

For this model set $f(i)=1_{i=0}$, for which the corresponding model $\Q_{n,\beta,\theta f}$ has sufficient statistic $h_0(G)$, the number of isolated vertices. This can be viewed as a penalty term which prefers or dislikes  isolated vertices depending on whether $\theta>0$ or $\theta<0$.
  Since $f$ is bounded, the asymptotics of the log normalizing constant follows from Corollary \ref{gwd}.
For this particular choice of  $f$, a direct calculation reveals that $$e^{Z(u,\theta  f)}=e^u+e^\theta -1,\quad m(u,\theta f)=\frac{ue^{u}}{e^u+e^\theta -1}.$$ 
Computation of $J(\beta,\theta f)$ can then be carried out numerically in a straightforward manner.

\subsection{Polynomial decay model}

For this model set $f(i)=i^\alpha$ for some known $\alpha \in [0,1]$. In this case the decay is at most linear by assumption, and so Corollary \ref{gwd} applies. Proceeding to compute $Z(u,\theta f)$ we have
$$e^{Z(u,\theta f)}=\sum_{i=0}^\infty \frac{u^i}{i!} e^{\theta i^\alpha}=\sum_{i=0}^\infty \frac{u^i}{i!}\sum_{j=0}^\infty \frac{\theta^j i^{\alpha j}}{j!}=e^{u+\theta}\sum_{i,j=0}^\infty \frac{e^{-u}u^i}{i!}\frac{e^{-\theta}\theta^j}{j!} i^{\alpha j}=e^{u+\theta}\E X^{\alpha Y},$$
where $X\sim Poisson (u)$ and $Y\sim Poisson (\theta)$ are mutually independent. A similar computation gives
$$m(u,\theta f)=\frac{\sum_{i=1}^\infty \frac{u^i}{(i-1)!} e^{\theta i^\alpha}}{\sum_{i=0}^\infty \frac{u^i}{i!} e^{\theta i^\alpha}}=u\frac{\sum_{i=0}^\infty \frac{u^i}{i!} e^{\theta (i+1)^\alpha}}{\sum_{i=0}^\infty \frac{u^i}{i!} e^{\theta i^\alpha}}=u\frac{e^{u+\theta} \E (X+1)^{\alpha Y}}{e^{u+\theta}\E X^{\alpha Y}}=u\frac{\E (X+1)^{\alpha Y}}{\E X^{\alpha Y}}.$$

Further simplification is not possible in general, and one has to use numerical methods to compute both $Z(u,\theta f)$ and $m(u,\theta f)$.

\section{Proofs of main results}\label{ldp}

 The main tool for proving our results is a large deviation principle for the empirical degree distribution $\mu_n^{G}$. To see how large deviation comes into the picture, note that the log normalizing constant of the model $\Q_{n,\beta, f}$ can be written as 
$$Z_n(\beta, f)=\log \E_{\P_{n,\beta}}e^{n\mu_n^{G}[f]},$$
where $\P_{n,\beta}$ is the Erd\"os-Renyi model with parameter $(\beta/n)$. By Varadhan's Lemma, this equates the problem to studying the large deviation of $\mu_n^{G}$ under the Erd\"os-Renyi $(\beta/n)$ model. A large deviation for the whole graph $G$ with respect to local weak convergence has recently been derived in \cite[Theorem 1.8]{BC}), which in particular gives a large deviation principle for $\mu_n^G$ with respect to the weak topology, as pointed out in \cite[Corollary 1.9]{BC}. The same large deviation was also obtained in \cite[Corollary 2.2]{DM} while studying large deviation for colored random graphs. 

The following definition introduces the rate function for this large deviation principle.
\begin{defn}\label{rate}

Let $\mathcal{S}\subset \P(\N_0)$ denote the set of all probability measures $\mu$ such that $\bar{\mu}=\sum_{i=1}^\infty i\mu(i)<\infty$.   Set the function $I_\beta:\P(\N_0)\mapsto [0,\infty]$ to be $+\infty$ if  $\mu\notin \cS$, and for $\mu\in\cS$ set 
\begin{align*}
I_\beta(\mu):=&\sum\limits_{i=0}^\infty \mu(i)\log (i!\mu(i))-\frac{\overline{\mu}}{2}\log (\overline{\mu}\beta)+\frac{\overline{\mu}+\beta}{2}\\
=&D(\mu||p_\beta)+\frac{1}{2}(\overline{\mu}-\beta)+\frac{\overline{\mu}}{2}\log \beta-\frac{\overline{\mu}}{2}\log \overline{\mu}
\end{align*}
 where $D(.||.)$ is the Kullback Leibler divergence, and $p_\beta$ is the Poisson distribution with parameter $\beta$.

\end{defn}

The following large deviation follows from \cite{BC,DM}:

\begin{thm}\label{done}

If $G$ is an Erd\"os-Renyi random graph with parameter $\beta/n$, then $\mu_n^G$ satisfies a large deviation principle on $\P(\N_0)$ with respect to weak topology, with speed $n$ and the good rate function $I_\beta(.)$.

\end{thm}
A direct application of the above large deviations result can be used to prove that
$$\lim_{n\rightarrow\infty}\frac{1}{n}Z_n(\beta, f)=\sup_{\mu\in \mathcal{S}}\{\mu[f]-I_\beta(\mu)\}$$
 when $f$ is a bounded function.   We now state three lemmas which will be used to extend this to all functions satisfying the conditions of Theorem \ref{thm:infty}.

\begin{lem}\label{technical}
For any function $f:\mathbb{N}_0\mapsto \R$ satisying \eqref{eq:consistency2} and any set $B\subset \P(\N_0)$ one has
$$\lim_{n\rightarrow\infty}\frac{1}{n}\log\E_{\P_{n,\beta}}e^{\sum_{i=0}^{n-1}h_i(G)f(i)}1\{\mu_n(G)\in B\}\le \sup_{\mu\in B\cap \cS}\{\mu[f]-I_\beta(\mu)\}.$$

\end{lem}

\begin{lem}\label{lem:technical2}
For finite positive real $\alpha$ and $f:\N_0\mapsto\R$ satisfying \eqref{eq:consistency2} we have
$$\sup_{\mu:I_\beta(\mu)-\mu[f]\le \alpha}\sum_{i=0}^\infty i\log i\mu(i)\le C,$$
where $C=C(\alpha,f,\beta)$ is a finite positive constant.
\end{lem}

\begin{lem}\label{above}
Let $f:\N_0\mapsto \R$ satisfy \eqref{eq:consistency2}. 

\begin{enumerate}
\item[(a)]
We have $$\sup_{\mu\in \mathcal{S}}\{\mu[f]-I_\beta(\mu)\}=J(\beta,f),$$ where $J({\beta},f)$ is as defined in \eqref{eq:j}. 
The supremum in this  definition is finite, and is attained over a finite set  of positive reals $\{u_1,\cdots,u_k\}$. Further,  any optimizing   $u$ satisfies the relation $u=\sqrt{\beta \overline{\sigma}_{u,f}}$.

\item[(b)]
For any $\varepsilon>0$ and $\psi$ satisfying \eqref{consistency} we have 
$$\sup_{\mu\in U^c}\{ \mu[f]-I_\beta(\mu)\}<\sup_{\mu \in \mathcal{S}}\{\mu[f]-I_\beta(\mu)\},$$
where $U:=\{\mu\in \P(\N_0):\min_{1\le l\le k}|\mu(\psi)-\sigma_{u_l}(\psi)|<\varepsilon\}.$

\end{enumerate}
\end{lem}

\subsection{Proofs of Theorem \ref{thm:infty}, Corollary \ref{gwd}, and Theorem \ref{gwd2}}\label{sec:proofs}

We now complete the proof of the main results of this paper, deferring the proof of the lemmas stated above to Section \ref{sub:proofs}.

\begin{proof}[Proof of Theorem \ref{thm:infty}]
\begin{enumerate}
\item[(a)]
To begin note that
\begin{align}\label{eq:free_upper}
e^{Z_n(\beta, f)}=\E_{\P_{n,\beta}} e^{ \sum_{i=0}^{n-1}h_i(G)f(i)},
\end{align}
which on taking $\log$, dividing by $n$ and letting $n\rightarrow\infty$ along with Lemma \ref{technical} gives
$$\limsup_{n\rightarrow\infty}\frac{1}{n}Z_n(\beta,f)\le \sup_{\mu\in \cS}\{\mu[f]-I_\beta(\mu)\},$$
and so we have verified the upper bound. The proof of the lower bound is split into two cases, depending on whether we are in case (i) or case (ii).

\begin{enumerate}
\item[(i)]
Define a function $T:\P(\N_0)\mapsto \R$ by 
\begin{align*}
T(\mu)=&\mu[f]\text{ if }I_\beta(\mu)<\infty,\\
=&0\text{ otherwise},
\end{align*}
 and use  Lemma \ref{lem:technical2} to note that $I_\beta(\mu)<\infty$ implies
 $\sum_{i=0}^\infty i\log i\mu(i)<\infty$. Also, since $f$ satisfies \eqref{consistency}, there exists $C_0<\infty$ such that $|f(i)|\le C_0i\log i$ for all $i\ge 0$. This immediately gives
 $$|T(\mu)|=|\mu(f)|\le C_0\sum_{i=0}^\infty i\log i\mu(i)<\infty.$$ 
Also for every $m\ge 1$ define the function $T_m:\P(\N_0)\mapsto \R$ by setting $T_m(\mu)= \sum_{i=0}^m f(i)\mu(i)$, and note that $T_m$ is continuous with respect to weak topology. We claim that for every positive real $\alpha$ and 
 $\delta>0$ we have
\begin{align}\label{eq:claim}
\lim_{m\rightarrow\infty}\sup_{\mu:I_\beta(\mu)\le \alpha}|T_m(\mu)-T(\mu)|=0.
\end{align}
To see this, fixing $\delta>0$ and invoking \eqref{consistency} we have $|f(i)|\le \delta i\log i$ for all $i> M(\delta)$. Thus for all $m\ge M(\delta)$ we have
$$|T_m(\mu)-T(\mu)|=|\sum_{i=m+1}^\infty f(i)\mu(i)|\le \delta \sum_{i=M(\delta)+1}^\infty i\log i\mu(i)\le \delta \sum_{i=0}^\infty i\log i\mu(i)\le \delta C(\alpha,0,\beta),$$
where the existence of $C(\alpha,0,\beta)$ follows from invoking Lemma \ref{lem:technical2} with $f\equiv 0$. Since $\delta>0$ is arbitrary, this verifies \eqref{eq:claim}.
\\

 We further claim that for every $\delta>0$ we have
\begin{align}\label{eq:claim2}
\lim_{m\rightarrow\infty}\limsup_{n\rightarrow\infty}\frac{1}{n}\log \P_{n,\beta}(|T_m(\mu_n^G)-T(\mu_n^G)|>\delta)=-\infty.
\end{align}
Indeed, with $\psi(i):=\sqrt{i\log i | f(i)|}$ we have $$|f(i)|\ll \psi(i)\ll i\log i,$$
and so there exists $M(\delta)$ such that for all $i\ge M(\delta)$ we have $|f(i)|\le \delta^2 \psi(i) $. Thus for $m\ge M(\delta)$ we have
$$|T_m(\mu)-T(\mu)|\le \delta^2 \sum_{i=M(\delta)+1}\psi(i)\mu(i)\le \delta^2\mu[\psi],$$
and so Markov's inequality gives
\begin{align*}
\P_{n,\beta}(|T_m(\mu_n^G)-T(\mu_n^G)|>\delta)\le &\P(\mu_n^G[\psi]>\frac{1}{\delta})\le e^{-\frac{n}{\delta}}\E_{\P_{n,\beta}} e^{n\mu_n^G[\psi]}.
\end{align*}
On taking $\log$, dividing by $n$ and letting $n\rightarrow\infty$ along with Lemma \ref{technical} we get
$$\limsup_{n\rightarrow\infty}\frac{1}{n}\log \P_{n,\beta}(|T_m(\mu_n^G)-T(\mu_n^G)|>\delta)\le -\frac{1}{\delta}+J(\psi,\beta).$$
Since $J(\psi,\beta)$ is finite and $\delta$ is arbitrary, \eqref{eq:claim2} follows.
\\

Given \eqref{eq:claim} and \eqref{eq:claim2}, it follows by \cite[Theorem 4.2.23]{DZ} and Theorem \ref{done} that $T(\mu_n^G)=\mu_n^G[f]$ satisfies a large deviation principle on $\R$ with the good rate function
$$\widetilde{I}(x):=\inf_{\mu\in \P(\N_0):\mu[f]=x}I_\beta(\mu).$$ Also Lemma \ref{technical}  gives
$$\frac{1}{n}\log \E_{\P_{n,\beta}} e^{2n \sum_{i=0}^{n-1}h_i(G)f(i)}=\frac{1}{n}\log \E_{\P_{n,\beta}} e^{2nT(\mu_n^G)}\le \sup_{\mu\in \cS}\{2\mu[f]-I_\beta(\mu)\}=J(\beta,2 f).$$
The right hand side above is finite by part (a) of Lemma \ref{above}, as $ 2f(.)$ satisfies \eqref{consistency}. This verifies \cite[(4.3.3.)]{DZ} with $\gamma=2$, and so by \cite[Theorem 4.3.1]{DZ} with $\phi(x)=x$ we have
\begin{align*}
\frac{1}{n}\log \E_{\P_{n,\beta}} e^{n\sum_{i=0}^{n-1}h_i(G)f(i)}=\frac{1}{n}\log \E_{\P_{n,\beta}} e^{n T(\mu_n^G)}
=& \sup_{x\in \R}\{\theta x-\widetilde{I}(x)\}\\=&\sup_{\mu\in \cS}\{ \mu(f)-I_\beta(\mu)\}=J(\beta, f),
\end{align*}
where the last equality again uses part (a) of Lemma \ref{above}. This completes the proof of part (a).

\item[(ii)]

Fixing $m\in \mathbb{N}$ one has

\begin{align*}
\mathbb{E}_{\P_{n,\beta}}e^{n\mu_n^{G}[f]}\ge & \mathbb{E}_{\P_{n,\beta}} e^{n\mu_n^{G}[f]}1_{\max_{1\le j\le n}d_j(G)\le m}\\
=   &\E_{\P_{n,\beta}}e^{\sum\limits_{i=0}^m h_i(G) f(i)}1_{\max_{1\le j\le n}d_j(G)\le m}\\
\ge &\E_{\P_{n,\beta}}e^{\sum\limits_{i=0}^m h_i(G) f(i)+f(m)\sum_{i=m+1}^{n-1}h_i(G)}1_{\max_{1\le j\le n}d_j(G)\le m},
\end{align*}
where the last inequality uses the fact that $f(i)\le f(0)=0$ for all $i$, as $f$ is non-increasing. This, along with an application of FKG inequality \cite[Prop. 1]{FKG} gives
$$\mathbb{E}_{\P_{n,\beta}}e^{n\mu_n^{G}[f]}\ge \E_{\P_{n,\beta}}e^{\sum_{i=0}^mh_i(G)f(i)+f(m)\sum_{i=m+1}^{n-1}h_i(G)}\P_{n,\beta}(d_1\le m)^n,$$
where we use the fact that the function $$G\mapsto \sum\limits_{i=0}^{n-1} h_i(G)\tilde{f}_m(i),\quad \tilde{f}_m(i)=\max(f(i),f(m)))$$ is  non-increasing on the space of graphs $\mathcal{G}_n$, as $\tilde{f}_m$ is non-increasing. Since $f(0)\le \tilde{f}_m(i)\le f(m)$, it follows that $\tilde{f}_m$ is bounded, and so an application of part (i) gives
\begin{align*}
\liminf\limits_{n\rightarrow\infty}\frac{1}{n}\log \E_{\P_{n,\beta}}e^{n\mu_n^{G}[\tilde{f}_m]}\ge& \sup_{\mu\in \mathcal{S}}\{\sum\limits_{i=0}^\infty \mu(i) \tilde{f}_m(i)-I_\beta(\mu)\}+\log p_{\beta}[0,m]\\
\ge & \sup_{\mu\in \mathcal{S}}\{\sum\limits_{i=0}^\infty \mu(i) f(i)-I_\beta(\mu)\}+\log p_{\beta}[0,m],
\end{align*}
where the last inequality uses the fact that $\tilde{f}_m\ge f$, and $p_\beta[0,m]$ is the probability that  a Poisson random variable with parameter $\beta$ is at most $m$.  
The lower bound follows on 
letting $m\rightarrow\infty$ and noting that $p_\beta[0,m]\rightarrow 1$.  Combining the upper and lower bound gives
$$\lim\limits_{n\rightarrow\infty}\frac{1}{n}Z_n(\beta,f)=\sup_{\mu\in \mathcal{S}}\{ \mu[f]-I(\mu)\}.$$

\end{enumerate}

\item[(b)]
By  part (a)  of Lemma \ref{above},  the supremum in the right hand side above is finite and equals $J(\beta, f)$ of \eqref{eq:j}, and  the set of optimizing $u$ in the definition of $J(\beta, f)$ has finite cardinality. Denoting this set by $\{u_1,u_2,\cdots,u_k\}$, let $$U:=\{\mu\in \P(\N_0)	:\min_{i=1}^k|\mu[\psi]-\sigma_{u_i, f}[\psi]|<\varepsilon\},$$
where $\varepsilon>0$ is fixed. 
Thus we have
\begin{align*}
&\limsup\limits_{n\rightarrow\infty}\frac{1}{n}\log \Q_{n,\beta, f}(\mu_n^G\in U^c)\\
\le &\limsup\limits_{n\rightarrow\infty}\frac{1}{n}\log \E_{\P_{n,\beta}}e^{n\mu_n^{G}[f]}1_{\mu_n^{G}\in U^c}-\liminf\limits_{n\rightarrow\infty}\frac{1}{n}\log Z_n(\beta,f)\\
\le & \sup_{\mu \in U^c\cap \cS}\{\mu[f]-I_\beta(\mu)\}-\sup_{\mu \in \mathcal{S}}\{\mu[f]-I_\beta(\mu)\},
\end{align*}
 where the last line uses  Lemma \ref{technical} with $B=U^c$, and part (a).
The last quantity above is negative by part (b) of Lemma \ref{above}, and so the conclusion follows.

\end{enumerate}
\end{proof}

\begin{proof}[Proof of Corollary \ref{gwd}]
\begin{enumerate}
Part(a) follows trivially from part (a) of Theorem \ref{thm:infty} on noting that the function $\theta f(.)$ satisfies \eqref{consistency} for all $\theta\in \R$. The continuity of the limiting log partition function follows from the fact that limit of convex functions is convex.

For part (b), setting $m:=\frac{1}{4}\min_{i=1}^k\bar{\sigma}_{u_i,\theta f}, M:=\max_{i=1}^k\bar{\sigma}_{u_i,\theta f}$, the desired conclusion follows from part (b) of Theorem \ref{thm:infty}.

\end{enumerate}
\end{proof}

\begin{proof}[Proof of Theorem \ref{gwd2}]

\begin{enumerate}
\item[(a)]
Since $f$ is non-decreasing and $\theta<0$, the function $\theta f$ is non-increasing and non-positive, and so an application of Theorem \ref{thm:infty} proves part (a).

\item[(b)]
It suffices to show that $\Q_{n,\beta,\theta f}(G=K_n)$ converges to $1$, as the desired conclusion about the log normalizing constant immediately follows. 

To this effect, using \eqref{eq:dir} there exists $M>4$ such that $f(i)- f(i-1)\ge M\log i$ for all $i\ge k_n:=\lfloor n/2\rfloor$, for all $n$ large enough.
We now claim that for all $r\in [0,n-1]$ we have 
\begin{align}\label{eq:crucial}
f(n-1)-f(n-1-r)\ge \frac{1}{4}Mr\log n
\end{align}
Indeed, if $r\le k_n$, then we have
\begin{align}\label{eq:crucial1}
\notag f(n-1)-f(n-1-r)=& \sum_{i=1}^{r}\Big(f(n-i)-f(n-i-1)\Big)\\
\ge &M\sum_{i=1}^r\log (n-i)\ge Mr\log (n/2).
\end{align}
On the other hand if $r> k_n$, using the monotonicity of $f$ along with \eqref{eq:crucial} gives
\begin{align}\label{eq:crucial2}
f(n-1)-f(n-1-r)\ge f(n-1)-f(n-1-k_n)\ge Mk_n\log(n/2).
\end{align}
Combining \eqref{eq:crucial1} and \eqref{eq:crucial2}, \eqref{eq:crucial} follows.

Thus if $G\in\cG_n$ is a graph with degree sequence $(d_1(G),\cdots,d_n(G))$, then setting $r_j(G):=n-1-d_j(G)$ for $G\in\cG_n$ we have
$$\sum_{j=1}^nf(n-1)-\sum_{j=1}^nf(d_j(G))\ge \frac{M\log n}{4}\sum_{j=1}^nr_j(G)= \frac{M\log n}{2}\Big(\frac{n(n-1)}{2}-E(G)\Big),$$
which immediately gives
\begin{align*}
\frac{\Q_{n,\beta,\theta f}(G)}{\Q_{n,\beta,\theta f}(K_n)}\le \Big(\frac{n}{\beta}\Big)^{R(G)} e^{-\frac{MR\log n}{2}R(G)},
\end{align*}
where $R(G):={n\choose 2}-E(G)$ for $G\in \cG_n$. This on summing gives
\begin{align*}
\Q_{n,\beta,\theta f}(R(G)\ge 1)\le \sum_{r=1}^{n\choose 2}{{n\choose 2}\choose r}e^{-\frac{Mr\log n}{2}}
\le \sum_{r=1}^{n\choose 2}n^{2r} e^{-\frac{Mr\log n}{2}}\le \sum_{r=1}^\infty \Big(n^2 e^{-\frac{M\log n}{2}}\Big)^r,
\end{align*}
which converges to $0$ as $n\rightarrow\infty$, as $M>4$.

\end{enumerate}
\end{proof}

\subsection{Proofs of Lemmas \ref{technical}, \ref{lem:technical2} and \ref{above}}\label{sub:proofs}

\begin{proof}[Proof of Lemma \ref{technical}]

Let $\cH_n$ to be the set of all degree frequency vectors $h(G)=(h_0(G),\cdots,h_{n-1}(G))$ on $n$ vertices as the graph $G$ varies in $\cG_n$. Fixing $\delta>0$ arbitrary we have
\begin{align}
\notag&\E_{\P_{n,\beta}}e^{\sum_{i=0}^{n-1} h_i(G)f(i)}1\{\mu_n(G)\in B\}\\
=&\mathbb{E}_{\P_{n,\beta}}e^{\sum_{i=0}^{n-1} h_i(G)f(i)}1\{\mu_n(G)\in B, E(G)\le \delta a_n\}
+\mathbb{E}_{\P_{n,\beta}}e^{\sum_{i=0}^{n-1} h_i(G)f(i)}1\{\mu_n^G\in B,E(G)>\delta a_n\}
\label{eq:appen1},
\end{align}
By \eqref{consistency} there exists $N=N(\delta)$ such that  $f(i)\le \frac{\delta}{4} i\log i$ for all $i> N(\delta)$, and so with $M:=M(\delta)=\max_{0\le i\le N}f(i)$ and $ a_n:=\frac{n(n-1)}{2}$ the second term in the right hand side of \eqref{eq:appen1} can be bounded by
\begin{align*}
e^{nM+\frac{\delta}{4}n^2\log n}\P_{n,\beta}(E(G)>\delta a_n)\le e^{nM+\frac{\delta}{4}n^2\log n}\sum_{r=\delta a_n}^{a_n}{a_n\choose r} \Big(\frac{\beta}{n}\Big)^{r}
\le  e^{nM+\frac{\delta}{4}n^2\log n}\times a_n2^{a_n}  \Big(\frac{\beta}{n}\Big)^{\delta a_n},
\end{align*}
which on taking $\log$, dividing by $n$, and letting $n\rightarrow\infty$ gives $-\infty$, and so we can ignore this term. The first term on the right hand side of \eqref{eq:appen1} can be written as
\begin{align*}
=&\sum_{{\bf h}\in \cH_n}N_n({\bf h}) e^{\sum_{i=0}^{n-1}h_i(G)f(i)}\Big(\frac{\beta}{n}\Big)^{E(G)}\Big(1-\frac{\beta}{n}\Big)^{{n\choose2}-E(G)} 1\{\mu_n^G\in B, E(G)\le \delta a_n\},
\end{align*}
where $N_n({\bf h})$ is the number of labeled graphs in $\cG_n$ whose degree frequency vector is ${\bf h}$.
It follows from \cite{M} that
 $$N_n({\bf h})\le \frac{(2r)!}{r!2^r\prod\limits_{i=0}^{n-1}i!^{h_i}}\times \frac{n!}{\prod\limits_{i=0}^{n-1}h_i!},$$
 where the extra factor $\frac{n!}{\prod\limits_{i=0}^{n-1}h_i!}$ accounts for the fact that for labeled graphs any relabeling between vertices with the same degree needs to be taken into account. Thus one has the following bound on the first term of the right hand side of \eqref{eq:appen1}:
\begin{align}\label{technical:bound1}
\mathbb{E}_{\P_{n,\beta}}e^{\sum_{i=0}^{n-1} h_i(G)f(i)}1\{\mu_n(G)\in B, E(G)\le \delta a_n \}
\le\sum_{{\bf h}\in \cH_n}\overline{N}_n( {\bf h)}1\{\mu_n(G)\in B, E(G)\le \delta a_n \},
\end{align}
where \begin{align*}
\overline{N}_n({\bf h}):=e^{\sum_{i=0}^{n-1}h_if(i)}\Big(\frac{\beta}{n}\Big)^{r}\Big(1-\frac{\beta}{n}\Big)^{{n\choose2}-r}\frac{(2r)!}{r!2^r\prod\limits_{i=0}^{n-1}i!^{h_i}}\times \frac{n!}{\prod\limits_{i=0}^{n-1}h_i!}
\end{align*}
 with $  r=\sum_{i=0}^{n-1}ih_i$. Using  Stirling's approximation one has
 \begin{align*}
C_2e^{-n}n^{n+1/2} \le n!\le C_1e^{-n}n^{n+1/2}
 \end{align*}
 for all $n\ge 1$, for some positive constants $C_1,C_2$ free of $n$.
Using this,  a direct computation gives
\begin{align*}
\overline{N}_n({\bf h})\le &e^{n\Big(1+o_n(1)\Big)\{\mu_n^G[f]-I_\beta(\mu_n^G)\}},
\end{align*}
which along with 
 \eqref{technical:bound1} gives
\begin{align}\label{eq:appen2}
\sum_{{\bf h}\in \cH_n}\overline{N}_n( {\bf h)}1\{\mu_n(G)\in B, E(G)\le \delta a_n \}\le e^{n\Big(1+o_n(1)\Big)\sup_{\mu\in B\cap \cS}\{\mu[f]-I_\beta(\mu)\}} \sum_{r=0}^{\delta a_n}|{\bf h}\in \cH_n:\sum_{i=0}^{n-1}ih_i=r|.
\end{align}
 Letting $p(r)$ denote the number of un-ordered partitions of $r$, we have the upper bound
$$|{\bf h}\in \cH_n:\sum_{i=0}^{n-1}ih_i=r|\le p(2r).$$
This is because given any such degree frequency vector ${\bf h}$, the corresponding un-ordered degree sequence sums up to $2r$, and so one can get a partition of $2r$ by dropping the vertices with degree $0$.
 Since  $$\lim_{r\rightarrow\infty}\frac{1}{\sqrt{r}}\log p(r)=\pi\sqrt{\frac{2}{3}},$$ (for a proof of this classical result see \cite{E} or \cite{HR} ),   taking logs, dividing by $n$ and taking $n\rightarrow\infty$   gives 
$$\limsup\limits_{n\rightarrow\infty}\frac{1}{n}\log \Big(\sum_{r=0}^{\delta a_n}|{\bf h}\in \cH_n:\sum_{i=0}^{n-1}ih_i=r|\Big)\le\pi \sqrt{\frac{2\delta}{3}}.$$
 This, along with \eqref{eq:appen1} and \eqref{eq:appen2} gives
$$\limsup_{n\rightarrow\infty}\frac{1}{n}\log \E_{\P_{n,\beta}} e^{\sum_{i=0}^{n-1}h_i(G)f(i)}1\{\mu_n^G\in B\}\le \sup_{\mu\in B\cap \cS}\{\mu[f]-I(\mu)\}+\pi \sqrt{\frac{2\delta}{3}},$$
from which the desired conclusion follows since $\delta>0$ is arbitrary.

\end{proof}

\begin{proof}[Proof of Lemma \ref{lem:technical2}]
Since $\log( i!)=\sum\limits_{k=1}^i\log k\ge \int\limits_{x=0}^i\log x dx=i\log i-i,$ we have
\begin{align}\label{n1}
\sum\limits_{i=0}^\infty \log( i!)\mu(i)\ge \sum\limits_{i=0}^\infty i\log i \mu(i)-\overline{\mu}
\end{align}
 
Also define $\sigma\in\mathcal{S}$  by $\sigma(i):=2^{-(i+1)}$ for $i\in\N_0$ and note that
\begin{align}\label{n2}
\sum\limits_{i=0}^\infty \mu(i)\log \mu(i)=D(\mu||\sigma)+\sum\limits_{i=0}^\infty \mu(i)\log \sigma(i)
\ge -(\overline{\mu}+1)\log 2.
\end{align}
Finally by \eqref{consistency}  there exists $M<\infty$ such that $f(i)\le M+ \frac{1}{4}i\log i$ for all $i\ge 0$. This gives
\begin{align}\label{n3}
\sum_{i=0}^\infty i\log i\mu(i)-\mu[f]\ge -M+\frac{3}{4}\sum_{i=0}^\infty i\log i\mu(i)\ge -M+\frac{3}{4}\bar{\mu}\log \bar{\mu},
\end{align}
where the last step uses Jensen's inequality.
Combining \eqref{n1}, \eqref{n2} and \eqref{n3} gives 
\begin{align}
\notag I_\beta(\mu)-\mu[f]=&\sum_{i=0}^\infty \log(i!)\mu(i)-\mu[f]+\sum_{i=0}^\infty \mu(i)\log\mu(i)-\frac{\bar{\mu}}{2}\log(\bar{\mu}\beta)+\frac{\bar{\mu}+\beta}{2}\\
\notag\ge &\sum_{i=0}^\infty i\log i\mu(i)-\bar{\mu}-\mu[f]-(\bar{\mu}+1)\log 2-\frac{\bar{\mu}}{2}\log(\bar{\mu}\beta)+\frac{\bar{\mu}+\beta}{2}\\
\notag\ge &-M+\frac{3}{4}\bar{\mu}\log\bar{\mu}-\bar{\mu}-(\bar{\mu}+1)\log 2-\frac{\bar{\mu}}{2}\log(\bar{\mu}\beta)+\frac{\bar{\mu}+\beta}{2}\\
=& -M+\frac{1}{4}\bar{\mu}\log \bar{\mu}-\overline{\mu}\Big(\log 2+\frac{3+\log \beta}{2}\Big)+\frac{\beta}{2}-\log 2=\phi_1(\overline{\mu}),\label{n4}
\end{align}
where $$\phi_1(x):= -M+\frac{1}{4}x\log x-x\Big(\log 2+\frac{3+\log \beta}{2}\Big)+\frac{\beta}{2}-\log 2.$$ Since $\phi_1(x)$ is continuous and diverges to $\infty$ as $x\rightarrow\infty$, it follows that 
$\phi_1(\overline{\mu})\le \alpha$ implies $\overline{\mu}\le K(\alpha)$ for some $K(\alpha)<\infty$. Thus we have
\begin{align*}
\alpha\ge I_\beta(\mu)-\mu[f]=&\sum_{i=0}^\infty \log(i!)\mu(i)-\mu[f]+\sum_{i=0}^\infty \mu(i)\log\mu(i)-\frac{\bar{\mu}}{2}\log(\bar{\mu}\beta)+\frac{\bar{\mu}+\beta}{2}\\
\ge &\sum_{i=0}^\infty i\log i\mu(i)-\bar{\mu}-\mu[f]-(\bar{\mu}+1)\log 2-\frac{\bar{\mu}}{2}\log(\bar{\mu}\beta)+\frac{\bar{\mu}+\beta}{2}\\
\ge &-M+\frac{3}{4}\sum_{i=0}^\infty i\log i\mu(i)-\bar{\mu}-(\bar{\mu}+1)\log2-\frac{\bar{\mu}}{2}\log(\bar{\mu}\beta)+\frac{\bar{\mu}+\beta}{2}\\
=&\frac{3}{4}\sum_{i=0}^\infty i\log i\mu(i)-\phi_2(x),
\end{align*}
where $\phi_2(x):=M+x+(x+1)\log 2+\frac{x}{2}\log(x\beta)-\frac{x+\beta}{2}$. Thus we have
\begin{align*}
\frac{3}{4}\sup_{\mu:I_\beta(\mu)-\mu[f]\le \alpha}\sum_{i=0}^\infty i\log i\mu(i)\le \alpha+\sup_{0\le x\le K(\alpha)}\phi_2(x),
\end{align*}
from which the conclusion of the Lemma follows.

\end{proof}

\begin{proof}[Proof of Lemma \ref{above}]
\begin{enumerate}
\item[(a)]
It suffices to consider the minimization of $\mu\mapsto \{I_\beta(\mu)-\mu[f]\}$ over $\cS$. To this effect, first note that $$\alpha:=\inf_{\mu\in \cS }\{I_\beta(\delta_0)-\delta_0[f]\}+1<\infty.$$ Indeed, taking  $\mu=\delta_0$ gives $$I_\beta(\delta_0)-\delta_0[f]= \beta/2-f(0)<\infty.$$ Thus  it suffices to minimize $\mu\mapsto \{I_\beta(\mu)-\mu[f]\}$ over the set $B_\alpha:=\{\mu:I_\beta(\mu)-\mu[f]\le \alpha\}$. By Lemma \ref{lem:technical2} we have 
\begin{align}\label{n5}
\sup_{\mu\in B_\alpha}\sum_{i=0}^\infty i\log i\mu(i)\le C(\alpha)<\infty,
\end{align} and so by Markov's inequality the set $B_\alpha$ is  tight with respect to weak topology.
 Let $\{\nu_k\}_{k\ge 1}$ be a sequence of measures in $B_\alpha$ such that $$\lim_{k\rightarrow\infty}\{I_\beta(\nu_k)-\nu_k[f]\}=\inf_{\mu \in B_\alpha \cap U^c}\{I_\beta(\mu)-\mu[f]\}.$$  Then by tightness of $B_\alpha$, there exists a subsequence which converges weakly to $\nu$, say. Without loss of generality assume the original sequence $\{\nu_k\}_{k\ge 1}$ converges weakly to $\nu$. Since $f$ satisfies \eqref{consistency}, invoking uniform integrability implied by \eqref{n5} it follows that $\nu_k(f)$ converges to $\nu(f)$.  This, along with the observation that $I_\beta(.)$ is lower semi continuous gives 
\begin{align*}
\inf_{\mu \in \cS}\{I_\beta(\mu)-\mu[f]\}=\lim_{k\rightarrow\infty}\{I_\beta(\nu_k)-\nu_k[f]\}\ge \{I_\beta(\nu)-\nu[f]\},
\end{align*}
and so $\nu$ attains the infimum.  Let $A\subset \cS$ be the set of all probability measures where the infimum is attained. Then for any $\mu\in \mathcal{S}$  and $\nu\in A$, by convexity of $\mathcal{S}$ we have $(1-t)\nu+t\mu\in \mathcal{S}$ for any $t\in [0,1]$. Thus with $u:=\sqrt{\overline{\mu}\beta}$ we have
 \begin{align*}
 &\frac{\partial}{\partial t}\Big[ I_\beta((1-t)\nu+t\mu)-(1-t)\nu[f]-t\mu[f]\Big]_{t=0}\ge 0\\
 \Leftrightarrow &\sum\limits_{i=0}^\infty \Big(1+\log \nu(i)+\log i!-\frac{i}{2}(1+\log \overline{\nu})-\frac{i}{2}\log \beta+\frac{i}{2}-f(i)\Big)(\mu(i)-\nu(i)))\ge 0\\
  \Leftrightarrow &\sum\limits_{i=0}^\infty \Big(\log \nu(i)+\log i!-i\log u-f(i)\Big)(\mu(i)-\nu(i))\ge 0\\
  \Leftrightarrow&D(\nu||\sigma_{u,f})+D(\mu||\nu)\le D(\mu||\sigma_{u,f}).
 \end{align*}
where $\sigma_{u,f}$ is as defined in definition (\ref{dd}). Since this holds for all $\mu\in \mathcal{S}$, setting $\mu=\sigma_{u,f}$ gives $D(\sigma_{u,f}||\nu)=0$, and so $\nu=\sigma_{u,f}$. Thus  $A\subset \Omega_f$, and consequently 
$$\sup_{\mu\in \cS}\{\mu[f]-I_\beta(\mu)\}=\sup_{u\ge 0}\{\sigma_{u,f}[f]-I_\beta(\sigma_{u,f})\}=J(\beta,f),$$
where the last equality follows by a simple algebra. It also follows from the proof that any $\sigma_{u,f}\in A$ must satisfy $u=\sqrt{\beta \overline{\sigma_{u,f}}}$.

Finally, to solve the optimization $u\mapsto \phi_1(u):=\sigma_{u,f}[f]-I_\beta(\sigma_{u,f})$ over $u\ge 0$, differentiating with respect to $u$ gives
$$\phi'_1(u)=-m'(u,f)\log \frac{u}{\sqrt{m(u,f)\beta}}.$$
Also setting $\phi_2(u):=\sum_{i=0}^\infty \frac{e^{f(i)}}{i!}u^i=e^{Z(u,f)}$ we have $m(u,f)=u\frac{ \phi_2'(u)}{\phi_2(u)}$, which on differentiating with respect to $u$ gives 
$$m'(u,f)=\frac{\phi_2(u)\phi_2'(u)+u\phi_2(u)\phi_2''(u)-u\phi_2'(u)^2}{\phi_2(u)^2},$$
and so $$\lim_{u\rightarrow0}m'(u,f)=\lim_{u\rightarrow 0}\frac{m(u,f)}{u}=\frac{\phi_2'(0)}{\phi_2(0)}=e^{f(1)-f(0)}>0.$$
This gives $\lim_{u\rightarrow 0}\phi_1'(u)=+\infty$, and so $u=0$ is not a local maxima of $\phi_1(.)$. Also it follows from \eqref{n5} that  optimizing measure $\mu$ satisfies $$\bar{\mu}\log\bar{\mu}\le \sum_{i=0}^\infty i\log i\mu(i)\le C(\alpha),$$
and so $m(u,f)\le C'$ for some finite constant $C'$. This along with the relation $u^2=m(u,f)\beta$ implies any optimizing $u$ is at most $\sqrt{C'\beta}$. Thus denoting $\widetilde{A}$ denote the subset of all positive reals $u$ which are global maximizers of the function $u\mapsto \phi_1(u)$, it follows that the set $\widetilde{A}$ is compact. Since an analytic non constant function on a bounded domain cannot have infinitely many minimizers, the set $\widetilde{A}$ must have finite cardinality. This completes the proof of part (a).

%
%
%
%
%

\item[(b)]

If $\inf_{\mu\in U^c}\{I_\beta(\mu)-\mu[f]\}=\infty$ then there is nothing to show. Assuming that $$\alpha':=\inf_{\mu\in U^c}\{I_\beta(\mu)-\mu[f]\}+1<\infty,$$  it suffices to minimize $\mu\mapsto \{I_\beta(\mu)-\mu[f]\}$ over $B_{\alpha'}\cap U^c$. Letting $\{\nu_k\}_{k\ge 1}$ be a sequence of measures in $B_{\alpha'}\cap U^c$ such that 
 
 $$\lim_{k\rightarrow\infty}\{I_\beta(\nu_k)-\nu_k[f]\}=\inf_{\mu\in B_{\alpha'} \cap U^c}\{I_\beta(\mu)-\mu[f]\},$$
 by a similar tightness and uniform integrability argument as in part (a) it follows that there exists a measure $\nu\in \cS$ such that $\{\nu_k\}_{k\ge 1}$ converges to $\nu$ weakly, and

$$\lim_{k\rightarrow\infty}\nu_k(f)=\nu(f),\quad \lim_{k\rightarrow\infty}\nu_k[\psi]=\nu[\psi].$$
 Since $\nu_k\in U^c$  and $\nu_k(\psi)$ converges to $\nu(\psi)$, we have $\nu\in U^c$. Since $U$ contains all the global minimizers of $\mu\mapsto \{I_\beta(\mu)-\mu[f]\}$, we have

\begin{align*}\inf_{\mu\in U^c}\{I_\beta(\mu)-\mu[f]\}=&\lim_{k\rightarrow\infty}\{I_\beta(\nu_k)-\nu_k(f)\}\text{ [By choice of $\{\nu_k\}_{k\ge 1}$]}\\
\ge &I_\beta(\nu)-\nu[f]\text{ [ By lower semi continuity of $I_\beta(.)$]}\\
>&\inf_{\mu\in \cS}\{I_\beta(\mu)-\mu[f]\},
\end{align*}
where the last step uses the fact that $\nu\in U^c$ is not in a global minimizer of $\mu\mapsto \{I_\beta(\mu)-\mu[f]\}$.
 This completes the proof of part (b).

\end{enumerate}
\end{proof}

\subsection{Proof of Theorems \ref{showconsistency1} and \ref{showconsistency2}}\label{fit_the_model}

\begin{proof}[Proof of Theorem \ref{showconsistency1}]

Differentiating with respect to $\theta, \log u,c$ and eliminating $c$ gives the least square equations
\begin{eqnarray*}
\theta \sum_{i=0}^L\Big(f(i)-\bar{f}\Big)^2+\log u\sum_{i=0}^L\Big(i-\frac{L}{2}\Big)\Big(f(i)-\bar{f}\Big)&=&\sum_{i=0}^L (f(i)-\bar{f})\log\frac{i! h_i(G)}{n}\\
\theta \sum_{i=0}^L\Big(i-\frac{L}{2}\Big)\Big(f(i)-\bar{f}\Big)+\log u\sum_{i=0}^L\Big(i-\frac{L}{2}\Big)^2&=&\sum_{i=0}^L \Big(i-\frac{L}{2}\Big)\log\frac{i! h_i(G)}{n},
\end{eqnarray*}
where $\bar{f}:=\frac{1}{L+1}\sum_{i=0}^Lf(i)$.
Thus we have the following matrix equation for the least square estimates:
\begin{align}\label{eq:matrix}
(\hat{\theta}_n, \log \hat{u}_n)A=\Big[ \sum_{i=0}^L \Big(f(i)-\bar{f}\Big)\log\frac{i! h_i(G)}{n},\quad \sum_{i=0}^L \Big(i-\frac{L}{2}\Big)\log\frac{i! h_i(G)}{n}\Big],\end{align} where $A$ is a $2\times 2$ matrix defined by
\[A=: 
  \left[ {\begin{array}{cc}
   \sum_{i=0}^L\Big(f(i)-\bar{f}\Big)^2 & \sum_{i=0}^L \Big(i-\frac{L}{2}\Big)\Big(f(i)-\bar{f}\Big)\\
        \sum_{i=0}^L\Big(i-\frac{L}{2}\Big)\Big(f(i)-\bar{f}\Big) & \sum_{i=0}^L\Big(i-\frac{L}{2}\Big)^2 \      \end{array} } \right].\]

Now, by part (b) of Theorem \ref{thm:infty}  it follows that there exists a finite set $\{u_1,u_2,\cdots u_k\}$ with $u_l>0$ such that any limit point of the measure $\mu_n^{G}$ is of the form $\sigma_{u_l,\theta f}$ for some $l,1\le l\le k$. This implies that
there exists a random variable $U_n$ taking values in $\{u_1,u_2,\cdots,u_k\}$ such that for all $i\in [0,L]$ we have
\begin{align}\label{universalworks}
\frac{h_i(G)}{n}-\frac{U_n^i}{i!}e^{\theta f(i)-Z(U_n,\theta f)}=o_P(1).
\end{align}
Plugging this estimate, Slutsky's Theorem implies
\begin{align*}
\sum_{i=0}^L (f(i)-\bar{f})\log \frac{i!h_{i}(G)}{n}=&\theta\sum_{i=0}^L (f(i)-\bar{f})^2+\log U_n\sum_{i=0}^L i(f(i)-\bar{f})+o_P(1),\\
\sum_{i=0}^L \Big(i-\frac{L}{2}\Big)\log \frac{i!h_{i}(G)}{n}=&\theta\sum_{i=0}^L \Big(i-\frac{L}{2}\Big)f(i)+\log U_n\sum_{i=0}^L \Big(i-\frac{L}{2}\Big)i+o_P(1),
\end{align*}
which along with \eqref{eq:matrix} gives
\begin{align}\label{eq:matrix2}
(\hat{\theta}_n-\theta ,\log \hat{u}_n-\log U_n)A=o_P(1).
\end{align}
We now claim that the minimum eigenvalue $\lambda_{\min}(A)$ is not $0$. Deferring the proof of the claim, let us first complete the proof of the Theorem. Given this claim, \eqref{eq:matrix2} implies $$\Big\|\hat{\theta}_n-\theta,\log \hat{u}_n-\log U_n\Big\|_2 \le \frac{1}{\lambda_{\min}(A)}  \Big\|(\hat{\theta}_n-\theta,\log \hat{u}_n-\log U_n)A\Big\|_2=o_P(1),$$
thus proving that $\hat{\theta}_n$ is consistent for $\theta$, and  $\hat{u}_n=U_n+o_P(1)$.

Since part (b) of Theorem \ref{thm:infty} with $\psi(i)=i$ implies 
$$m(U_n,\theta f)-\frac{2E(G)}{n}=o_P(1),$$  we have
$$\frac{n\hat{u_n^2}}{2E(G)}=\frac{U_n^2}{m(U_n,\theta f)}+o_P(1)=\beta+o_P(1).$$
 where the last equality invokes the relation $U_n^2=\beta m(U_n,\theta f)$. This shows consistency of $\hat{\beta}_n$ for $\beta$ as well.

It thus remains to verify the claim that $\lambda_{\min}(A)$ is not $0$. To see this, note that if $\lambda_{\min}(A)=0$, then $|A|=0$, which gives
\begin{align*}
\sum_{i=0}^L(f(i)-\bar{f})^2\sum_{i=0}^L \Big(i-\frac{L}{2}\Big)^2=\Big[\sum_{i=0}^L \Big(i-\frac{L}{2}\Big)\Big(f(i)-\bar{f}\Big)\Big]^2.
\end{align*}
Thus equality holds in the Cauchy-Schwarz inequality, which implies $f(i)-\bar{f}=b\Big(i-\frac{L}{2}\Big)$ for some $b\in \R$. But then $f(0)=0$ forces $f(i)=bi$ for all $i\ge 0$, a contradiction. This completes the proof of the Theorem.
\end{proof}

\begin{proof}[Proof of Theorem \ref{showconsistency2}]

As before there exists a random variable $U_n$ taking values in a finite set $\{u_1,\cdots,u_k\}$ such that (\ref{universalworks}) holds, which gives
\begin{align*}
\frac{i!h_i(G)}{h_0(G)}-e^{f(i)}U_n^i=o_P(1).
\end{align*}
On taking $\log$  and using the definition of $\hat{f}_n(i)$ as in the theorem, this gives
$$\hat{f}_n(i)+i\log\hat{u}_n=\log \Big[\frac{i!h_i(G)}{h_0(G)}\Big]=f(i)+i\log U_n+o_P(1).$$
Thus to complete the proof it suffices to show that $\hat{u}_n-U_n=o_P(1)$. To prove this, first note that  part (b) of Theorem \ref{gwd2} gives
$$\frac{2E(G)}{n}-m(U_n,f)=o_P(1).$$
Since  one has $U_n^2=m(U_n,f)$ as well, it readily follows that $$\hat{u}_n=\sqrt{\frac{2E(G)}{n}}=\sqrt{m(U_n,f)}+o_P(1)=U_n+o_P(1),$$
thus completing the proof of the theorem.

\end{proof}

\section{Acknowledgement}
The content of this paper also appears in my Ph.D. thesis under the guidance of Prof. Persi Diaconis. I would also like to thank  Prof. Sourav Chatterjee and Prof. Amir Dembo for helpful discussions on this paper. This paper also benefitted from some helpful discussions during an AIM conference on ERGMs held in June 2013, for which I would like to thank all participants of the conference. The presentation of this paper has also benefitted from the recommendations of an anonymous reviewer, and the Associate Editor.

\end{document}